\documentclass[12pt]{article}
\usepackage{t1enc}
\usepackage[latin1]{inputenc}
\usepackage[english]{babel}
\usepackage{amssymb}
\usepackage{amsmath}

\usepackage{amsmath}
\usepackage{graphics}

\usepackage{amsthm, amssymb}
\usepackage{amsfonts}
\usepackage{epsfig,multicol}

\newtheorem{theorem}[subsection]{Theorem}
\newtheorem{proposition}[subsection]{Proposition}
\newtheorem{lemma}[subsection]{Lemma}
\newtheorem{corollary}[subsection]{Corollary}
\newtheorem{definition}[subsection]{Definition}

\newtheorem{remark}[subsubsection]{Remark}

\newfont{\gothic} { ygoth scaled \magstep{1.5}}

\def\<{\langle}
\def\>{\rangle}

\def\hpic #1 #2 {\mbox{$\begin{array}[c]{l} \epsfig{file=#1,height=#2}
\end{array}$}}
 
\def\vpic #1 #2 {\mbox{$\begin{array}[c]{l} \epsfig{file=#1,width=#2}
\end{array}$}}

\title{An orthogonal approach to  the subfactor of a planar algebra}
\author{Vaughan Jones\thanks{Supported in part by NSF Grant DMS0401734, Auckland University and the NZIMA}, Dimitri Shlyakhtenko and Kevin Walker}

\begin{document}
\DeclareGraphicsExtensions{.eps}

\maketitle 

\begin{abstract}
By changing to an orthogonal basis, we give a short proof that the subfactor of the graded algebra of a planar algebra  
reproduces the planar algebra. 
\end{abstract}
\section{Introduction}
Starting from a subfactor planar algebra, a construction was given in \cite{gjs:free} of a tower of II$_1$ factors
whose standard invariant is precisely the given planar algebra. The construction was entirely in
terms of planar diagrams and gave a diagrammatic reproof of a result of Popa in \cite{popa:standardlattice}. The inspiration
for the paper was from the theory of large random matrices where expected values of words on 
random matrices give rise to a trace (see \cite{voiculescu:symmetries}) on the algebra of noncommutative polynomials. 
Since that trace is definable entirely in terms of planar pictures it was easy to generalise it to an
arbitrary planar algebra, giving the planar algebra a concatenation multiplication to match that
of noncommutative polynomials. Unfortunately, though the algbera structure is very straightforward,
the inner product is not always easy to work with as words of different length are not orthogonal.
In this paper we use a simple diagrammatic orthogonalisation discovered by the third author to
reprove the II$_1$ factor results of \cite{gjs:free} in a direct and simple way without the use of full Fock space or
graph C$^*$-algebras. One may capitalise on the advantages of orthogonalisation since the multiplication
does not actually become much more complicated when transported to the orthogonal basis.
We present the results by beginning with the orthogonal picture and giving a complete proof of
the tower result. Then we show that this orthogonal structure is actually isomorphic to that of \cite{gjs:free}.
The same result   was obtained simultaneously and independently by Kodiyalam and Sunder in \cite{sunder}.
\section{Setup.} 
  Let $\mathfrak P = (P_{n})_{n=0,1,2,\cdots}$ be a subfactor planar algebra. Let $Gr_k(\mathfrak P)$ be the graded vector space
  $\oplus_{n\geq 0}P_{n+k}$ equipped with the prehilbert space inner product $<,>$ making 
  it an orthogonal direct sum and for which, within $P_{n,k}$, $<a,b>=\delta^{-k}\vpic {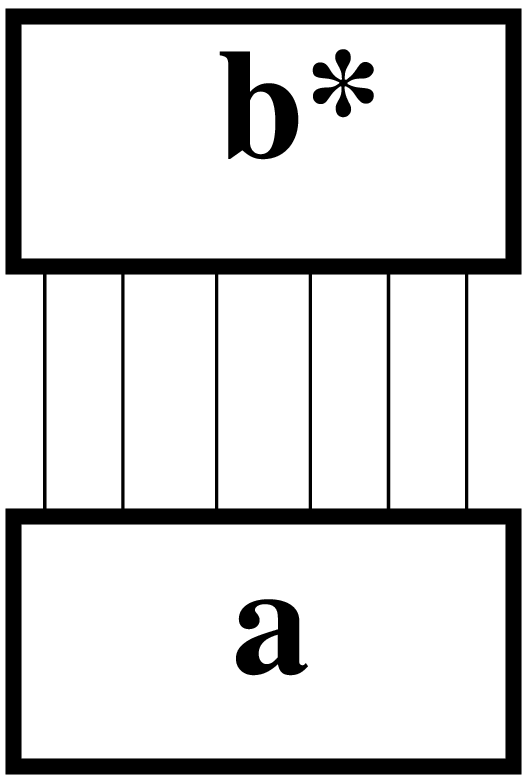} {0.5in} .$
  Write $P_{n,k}$ for $P_{n+k}$ when it is considered as the $n$-graded part of  $Gr_k(\mathfrak P)$.
  We will attempt to keep the pictures as uncluttered as possible by using several conventions and being
  as implicit as possible. Shadings for instance will always be implicit and we will eliminate the outside
  boundary disc whenever convenient. An element
  $a\in P_{n,k}$ will  be represented whenever possible in a picture as: $\vpic{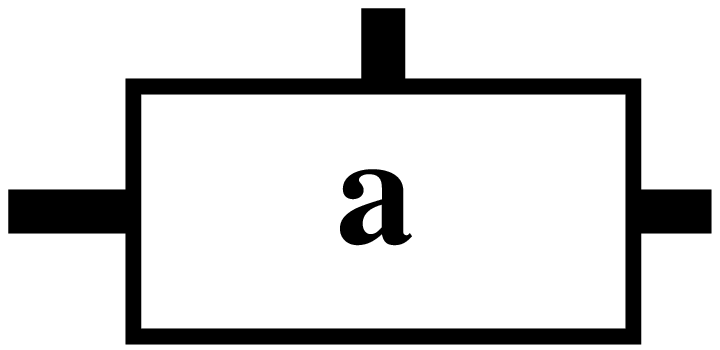} {1in} $ where
  the thick lines to the left and right of the box represent $k$ lines and the think line at the top
  represents $2n$ lines. If the multiple lines have to be divided into groups the number off lines in
  each group will be indicated to the minimal extent necessary. The distinguished first interval in
  a box will always be the top left of the box. Thus the inner product above of $<a,b>$ for elements
  of $P_{n,k}$ will be \vpic{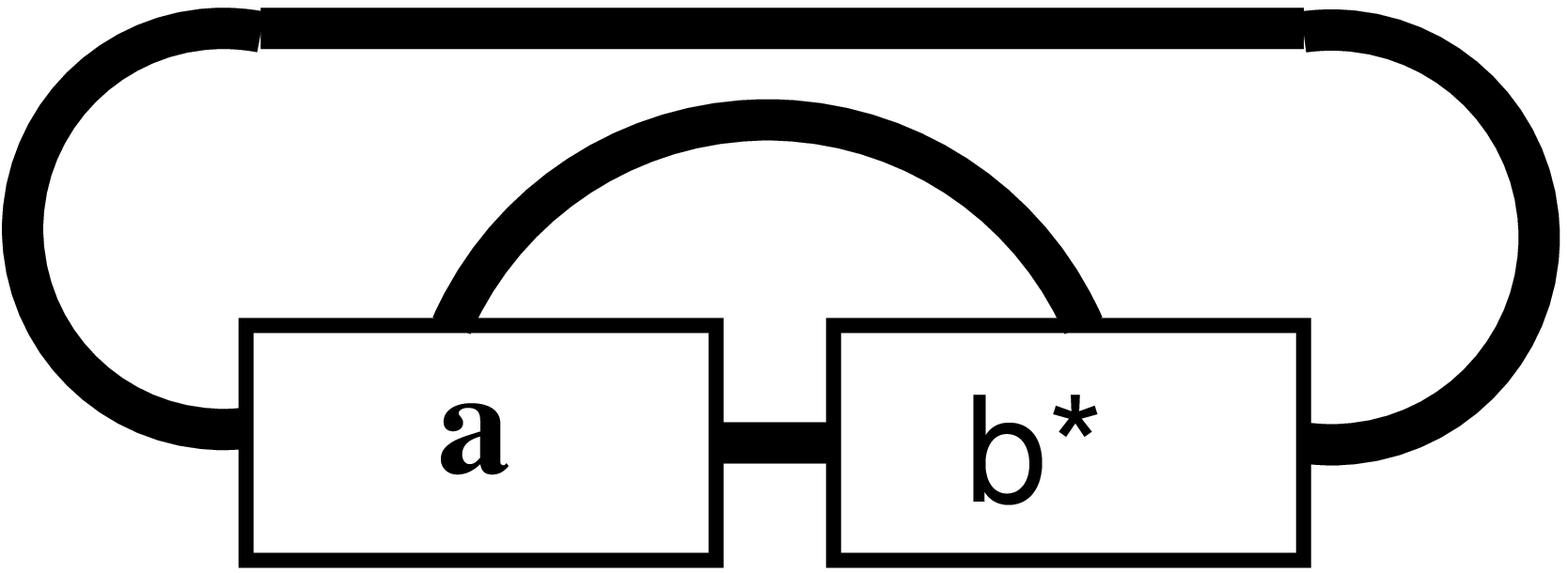} {1.2in} . In the original works on planar algebras (e.g. \cite{jones:planar}), 
  each $P_k$ ia an associative *-algebra whose product, with these conventions, would be to 
  view $P_k$ as $P_{0,k}$ and  $ab=\vpic {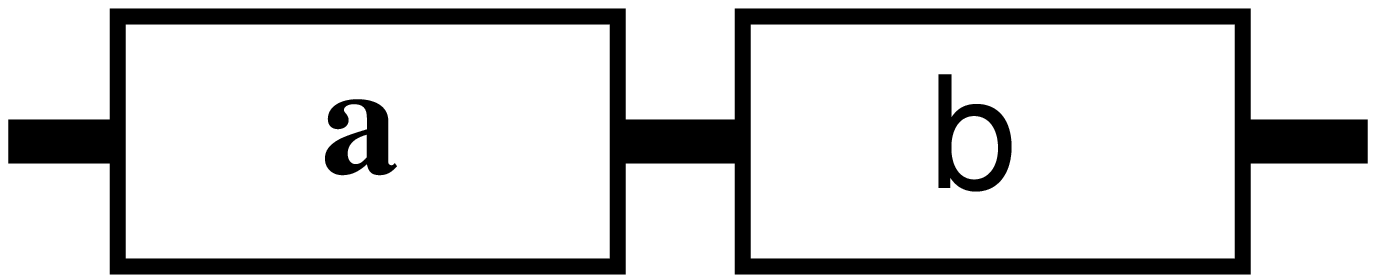} {1in} $. There are unital inclusions of $P_k$ in $P_{k+1}$
  by identifying $a$ with \vpic {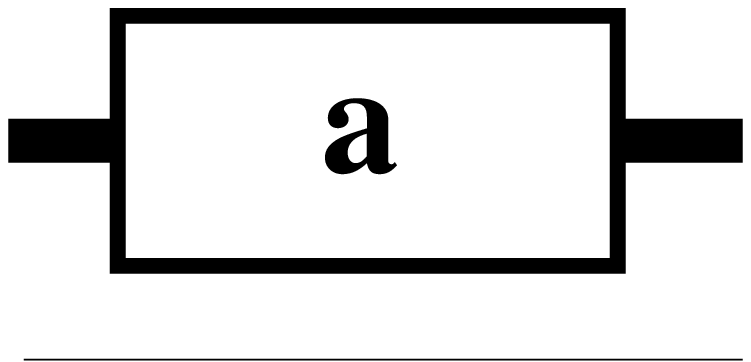} {1in} .  The identity element of $P_k$ is thus represented
  by a single thick horizontal line. It is also the identity element of $Gr_k(\mathfrak P)$.  The trace (often
  called the Markov trace) $tr$ on$P_k$ is normalised so as to be compatible
  with the inclusions by $tr(a)=\delta^{-k}<a,1>$. We extend this trace to $Gr_k(\mathfrak P)$ by the
  same formula so that the trace of an element is the Markov trace of its zero-graded piece. Each $P_n$
  is a finite dimensional C$^*$-algebra whose norms are also compatible with the inclusions. 
  
  \section{*-Algebra structure on $Gr_k(\mathfrak P)$}
  \begin{definition} If $a\in P_{m,k} $ and $b\in P_{n,k}$ are elements of $Gr_k(\mathfrak P)$ we
  define their product to be 
  $$ a\star b = \sum_{i=0}^{min(2m,2n)} \vpic {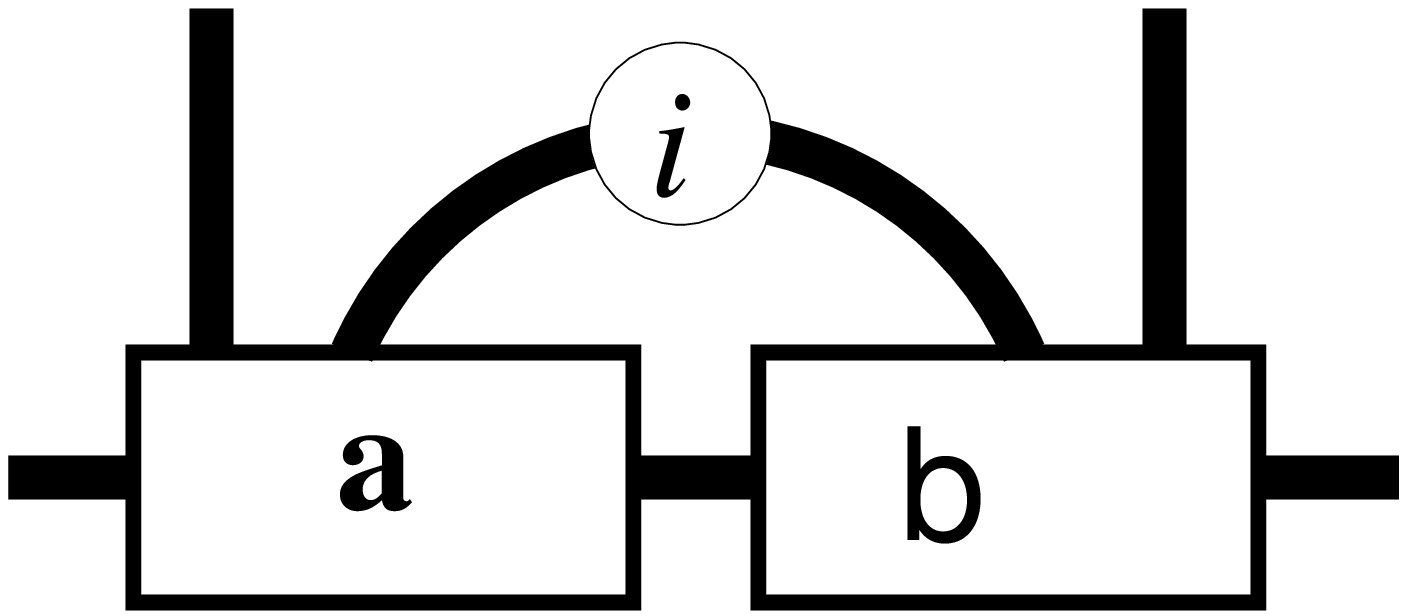} {1.5in} $$
  
  (where the $i$ means there are $i$ parallel strings. The numbers of other parallel strings 
  are then implicitly defined by our conventions.)
  
  The *-structure on $P_{n,k}$ is just the involution coming from the subfactor planar algebra.
  \end{definition}
  \begin{proposition} $(Gr_k(\mathfrak P), \star, ^*)$ is an associative *-algebra.
  \end{proposition}
  \begin{proof}
The property $(a\star b)^*=b^* \star a^*$ is immediate from the properties of a planar *-algebra.
For associativity note that both $a\star(b\star c)$ and $(a\star b)\star c$ are given by the sum over
all 
epi (see section \ref{basis}) diagrams where no strand has both of its endpoints of $a$,
or both of its endpoints on $b$, or both of its endpoints on $c$.
Here are two typical examples:\\
\vpic{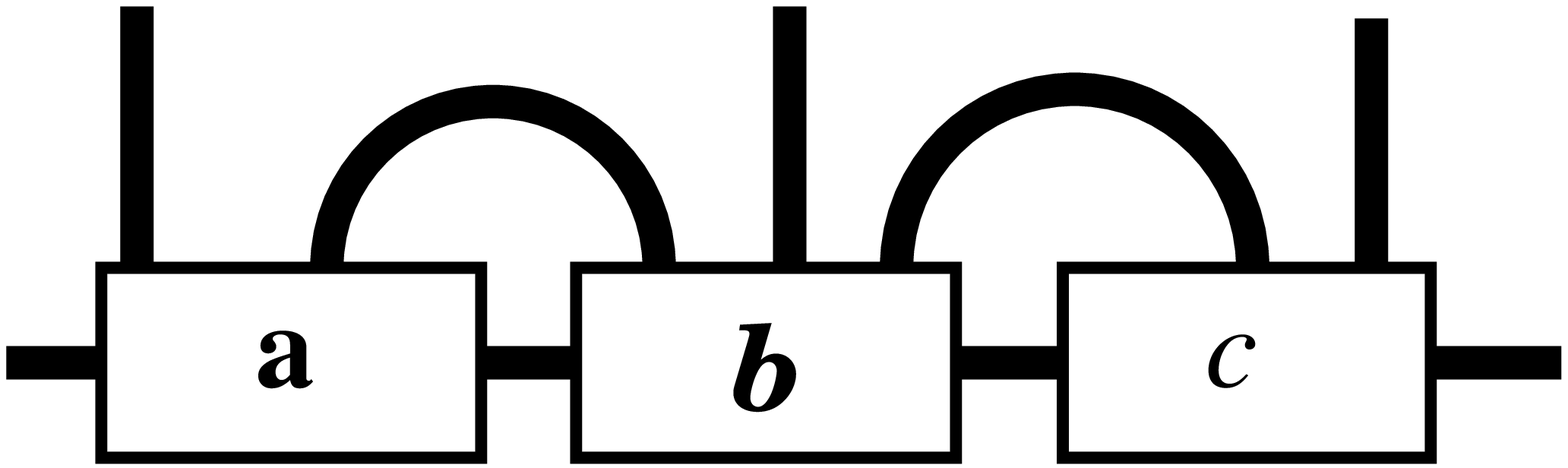} {1.9in}  \hspace{20pt} \vpic{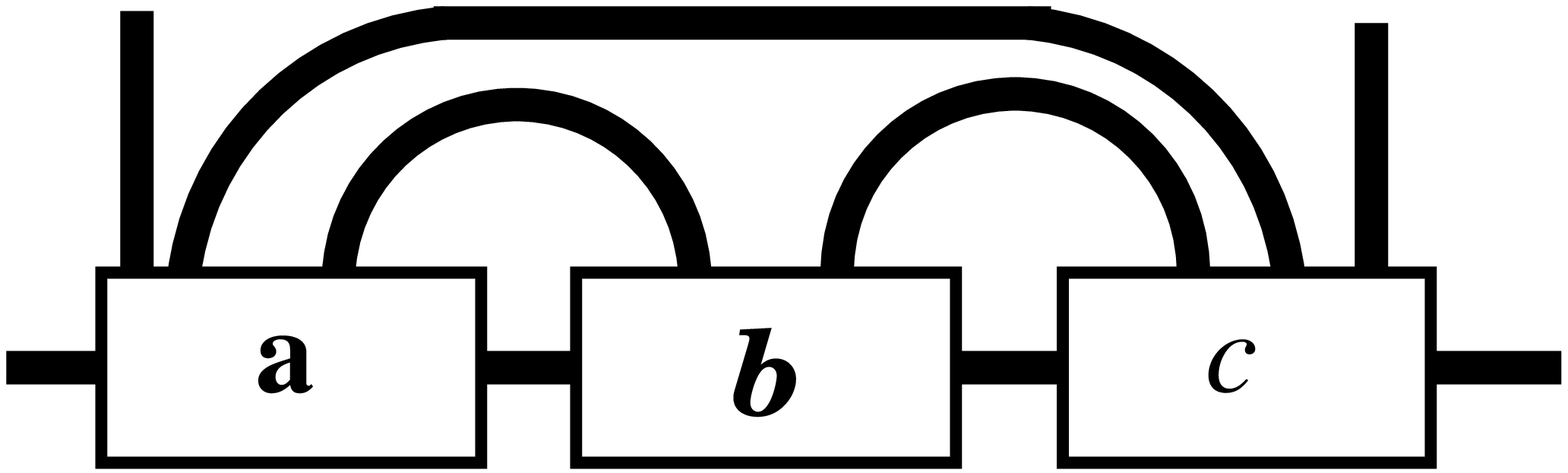} {1.9in}

  \end{proof}
  
  The inner product $<a,b>$ is clearly equal to $tr(ab^*)$ and is positive definite by 
  definition. We would like to perform the GNS
  construction but since there is no C$^*$-algebra available we need to show by hand
  that left (and hence right) multiplication by elements of $Gr_k(\mathfrak P)$ is bounded.
  
  \begin{theorem} Let $a\in Gr_k(\mathfrak P)$. Then the map $L_a:Gr_k(\mathfrak P)\rightarrow Gr_k(\mathfrak P)$,
 defined by  $L_a(\xi) =a\star \xi$ is bounded for
  the prehilbert space structure.
  \end{theorem}
  \begin{proof} We may suppose $a\in P_{n,k}$ for some $n$. Then $L_a$ is a sum of $2n+1$
  maps $L_a^i$ from
  an orthogonal direct sum of finite dimensional Hilbert spaces to another, respecting the
  orthogonal decomposition,$L_a^i$ being the map defined by the $i$th. term in the
  sum defining $\star$. Thus it suffices to show that the norm of the map $L_a^i: P_{m,k}
  \rightarrow P_{m+n-i,k}$ is bounded independently of $m$, the number of $i$ values being
  $\leq 2n+1$.  Clearly we may suppose that
  $m>>n+k$ which simplifies the number of pictures to be considered. So if $b\in P_{m,k}$
  we must estimate $<ab,ab>$ which is the following tangle:
  
 $\delta^{-k}$  \vpic{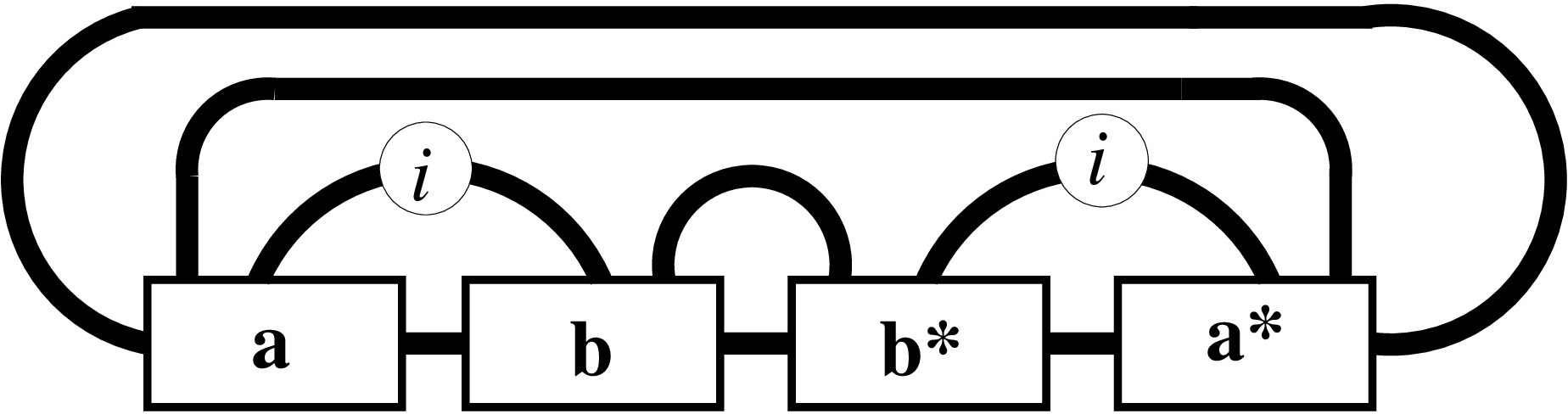} {3in}
  
  We may suppose $i\leq n$ since the norm of an operator is equal to that of its adjoint and the roles of $i$ and $2n-i$ are reversed in
  going between $L_a^i$ and $(L_a^i)*$.
  Then we may isotope the picture, putting $a$ and $b$ in boxes with the 
  same number ($k+n$ and $k+m$ respectively) of boundary points on the top and bottom, also 
  possibly rotating them, to obtain the following equivalent tangle:
  
$\delta^{-k}$  \vpic{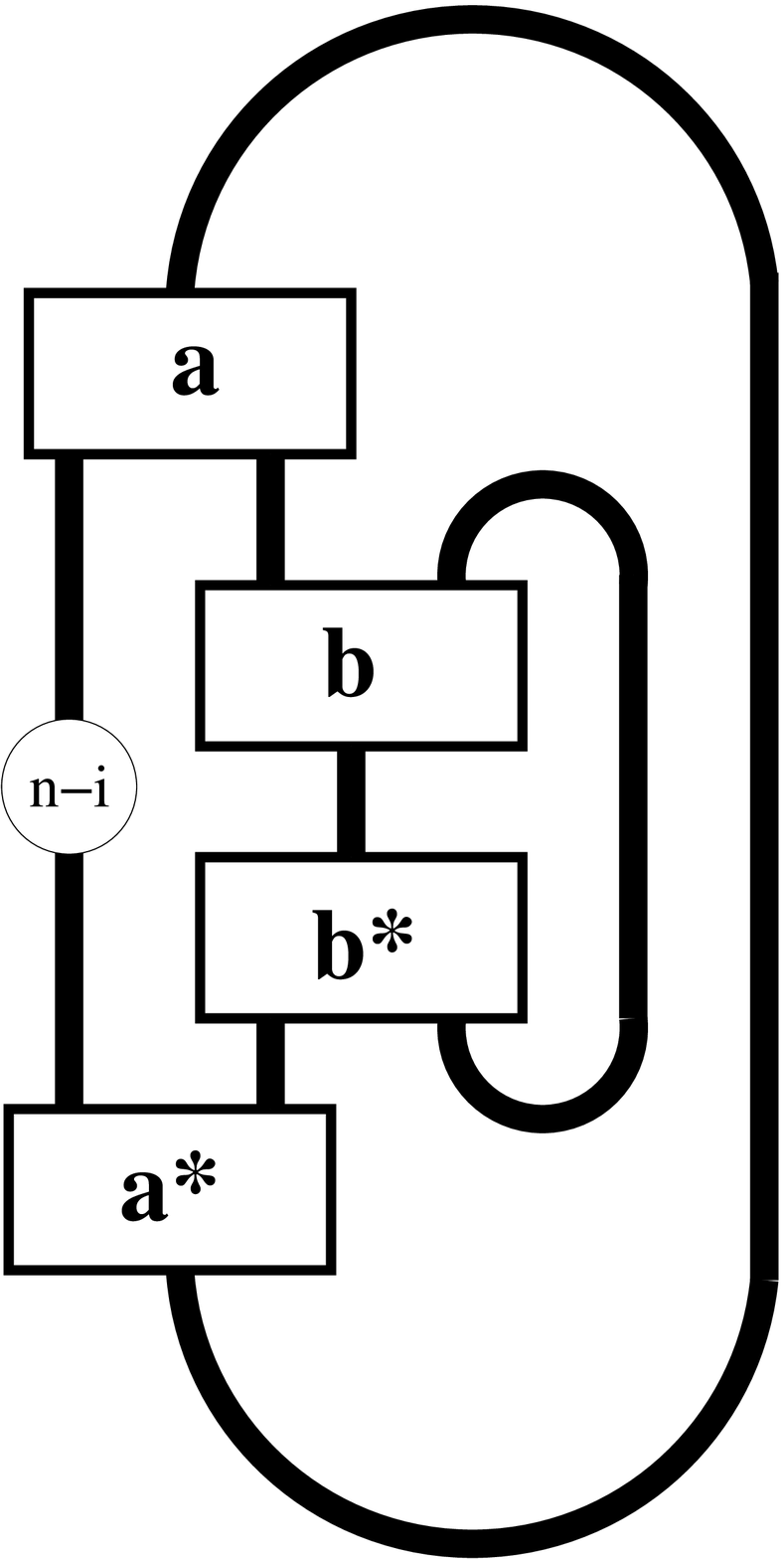} {1in}
  
  Note that the multiplicities of all the strings are determined by the $"n-i"$ and our
  conventions.
  
  Neglecting powers of $\delta$ which do not involve $m$ we see $<\tilde a \tilde b,\tilde a \tilde b>$ where $\tilde a$ is $a$ with $m-i$ strings to the right
  and $\tilde b$ is $b$ with $n-i$ strings to the left. The strings to the right do not change the norm
  of $a$ (as an element of the finite dimensional C$^*$-algebra
  $P_{2k+2n}$) by the uniqueness of the C$^*$-norm. The ($L^2$-) norm of 
  $\tilde b$ differs from that of $b$ by an $m$-independent power of $\delta$. Hence we are done.
  
  
  
  
  \end{proof}
  $Gr_k(\mathfrak P)$ has thus  been shown to be what is sometimes called a ``Hilbert Algebra''
  or ``unitary algebra''.
  
  \section{The von Neumann algebras $M_k$.}
  
  \begin{definition} Let $M_k$ be the finite von Neumann algebra on the Hilbert  space completion
  of $Gr(\mathfrak P_k)$ generated by left multiplication by the $L_a$.
  \end{definition}
  
  Since right multiplication is also bounded, the identity in $P_{0,k}$ is a cyclic and separating
  trace vector for $M_k$ defining the faithful trace $tr$ as usual, and the right multiplications generate
  the commutant of $M_k$. We shall first show that each $M_k$ is a factor. 
  
  \begin{definition} The element $\cup_k\in P_{2,k}$ will be the following:
  \vpic {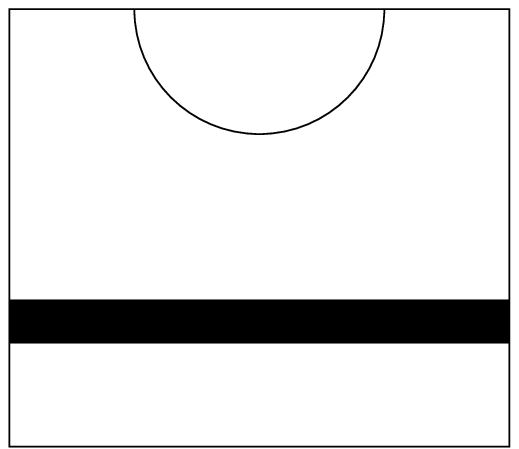} {0.5in} (where we have included the boundary to avoid disembodiment).
  The subalgebra of $Gr(\mathfrak P_k)$ generated by $\cup_k$ will be denoted $\mathfrak A_k$ and
  its weak closure in $M_k$ will be called $A_k$.
  \end{definition}
  
  \begin{definition} For $x\in P_{n,k}$, $n\geq 0$ and $p,q\geq 0$ let $x_{p,q}\in P_{n+p+q,k}$ be
  $\frac{1}{(\sqrt \delta)^{p+q}} $ times the following:
  
  \vpic{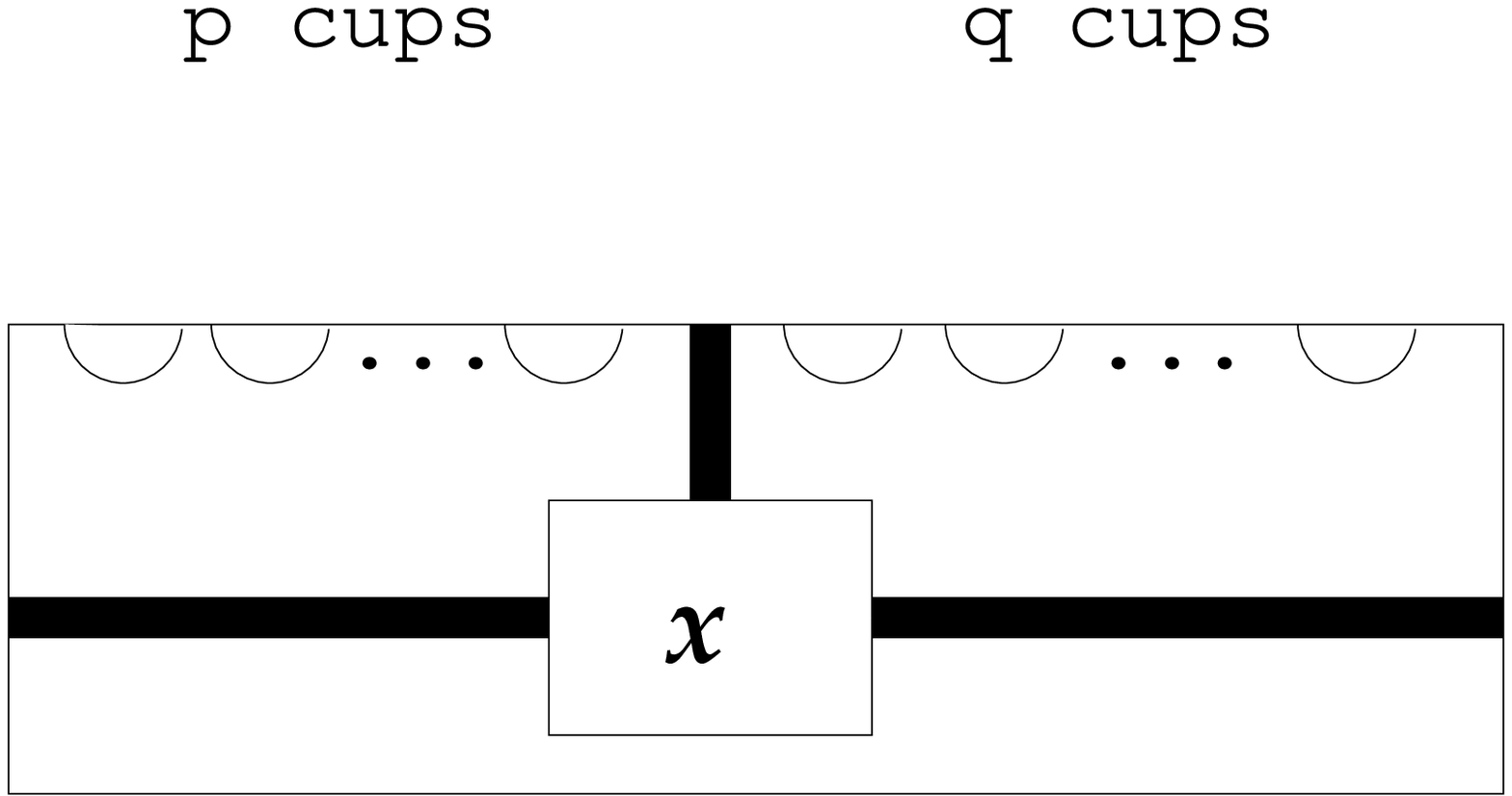} {3in}
   \end{definition}
  
  \begin{definition}Let $W_0 = \{0\} \subset P_{0,k}$ and for each $k\geq 1$ let $W_n$ be the
  span in $P_{n,k}$  of $\{x_{1,0}\}$ and $\{x_{0,1}\}$ for $x\in P_{n-1,k}$. For $n\geq 0$ let
  $V_n=W_n^\perp$.
  \end{definition}
  
  \begin{lemma} For $x\in P_{n,k}$, $n\geq 1$,  
  $$x\in V_n \iff  \vpic{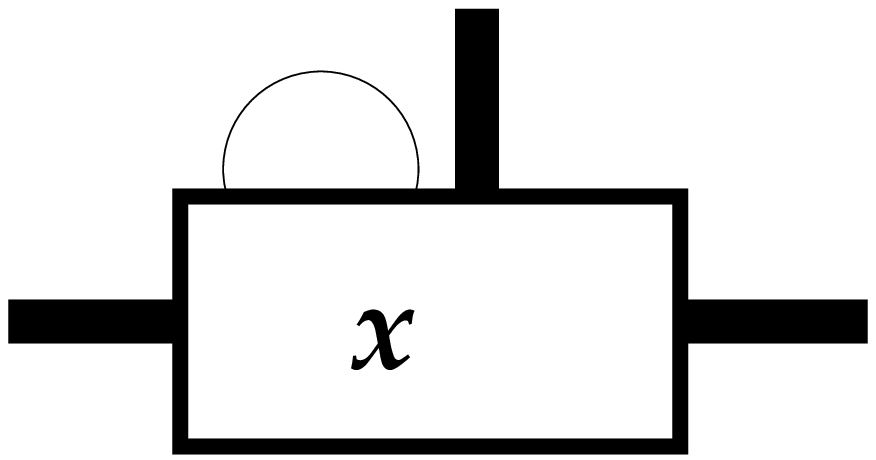} {1in} =0=\vpic{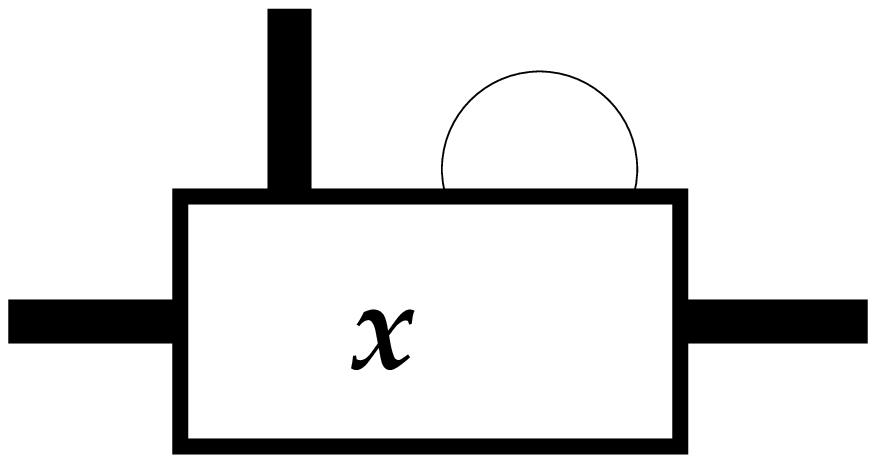} {1in}  .$$
  \end{lemma}
  \begin{proof}
  Taking the inner product of these two elements with an arbitrary element in $P_{n-1,k}$ we
  see the inner product of  element s in $W_{n,k}$ with $x$.
  \end{proof}
  \begin{corollary} Let $v\in V_m$, $m\geq 0$ and $v'\in V_n$, $n\geq 1$, then 
  $$<v_{p,q},v'_{p',q'}>=\left\{
  \begin{array} { ll}
    <v,v'> & \textrm{ if  } p=p',q=q' \\
    0         & \textrm{ otherwise}
    \end{array} \right. $$
    \end{corollary}
    \begin{proof} If either $p\neq p' $ or $q\neq q'$, the left or rightmost pair of boundary points
    of $v$ or $v'$ will be capped off to give zero.
    \end{proof}
    \begin{corollary} \label{formula} If $v\in V_n$ for $n>0$ is a unit vector, the $v_{p,q}$
    are an orthonormal basis for the $\mathfrak A-\mathfrak A$ bimodule $\mathfrak A v \mathfrak A$.
    \end{corollary}
    \begin{proof} By the previous lemma it suffices to show that the span of the $v_{p,q}$ is
    invariant under left and right multiplication (using $\star$) by $\cup_k$. In fact we have the
    formula:\\
    $$\cup_k \star v_{p,q}=\left\{
  \begin{array} { ll}
  \sqrt \delta v_{1,q} + v_{0,q} & \textrm{ if  } p=0 \\
   \sqrt \delta v_{p+1,q} + v_{p,q}+     \sqrt \delta v_{p-1,q}     & \textrm{ otherwise}
    \end{array} \right. $$
    And there is an obvious corresponding formula for right multiplication by $\cup_k$.
    \end{proof}
    \begin{lemma}
    The linear span of all the $v_{p,q}$ for $v\in V_n$ for all $n$ is $Gr(\mathfrak P_k)$.
    
    \end{lemma}
    \begin{proof} By a simple induction on $n$ these vectors span $P_{n,k}$ for all $k$.
    \end{proof}
    Let us summarise all we have learnt using the unilateral shift $S$ (with $S^*S=1$) on $\ell^2(\mathbb N)$.
    
    \begin{theorem}
    Suppose $\delta >1$. As an $A-A$ bimodule, 
    $$L^2(M_k)=P_{0,k}\otimes\ell^2(\mathbb N)\oplus \{\mathfrak H \otimes \ell^2(\mathbb N)\otimes \ell^2(\mathbb N)\}$$
    with $\cup_k$ acting on the left and right on $P_{0,k}\otimes \ell^2(\mathbb N)$ by $id\otimes(\sqrt \delta (S+S^*)+SS^*)$, on the left
    on $\mathfrak H\otimes  \ell^2(\mathbb N)\otimes \ell^2(\mathbb N)$ by 
    $id\otimes(\sqrt \delta (S+S^*) +1)\otimes id$ and on the right on
     $\mathfrak H\otimes  \ell^2(\mathbb N)\otimes \ell^2(\mathbb N)$ by $id\otimes id\otimes(\sqrt \delta (S+S^*) +1)$.
    ($\mathfrak H$ is an auxiliary infinite dimensional Hilbert space.)
    \end{theorem}
    
    \begin{proof} Obviously $P_{0,k}$ commutes with $A$ so the first term in the direct sum is the result of
    a simple calculation. Choosing an orthonormal basis for each $V_n$ gives the rest by corollary \ref{formula}.
    
    \end{proof}
    
    \begin{corollary}
    $A'\cap M_k=AP_{0,k}$
    \end{corollary}
    \begin{proof} It suffices to show that  no non-zero $\xi \in \ell^2(\mathbb N)\otimes \ell^2(\mathbb N)$ satisfies 
    $(S+S^*)\xi= \xi(S+S^*)$.  But such a $\xi$ would be a Hilbert-Schmidt operator on $\ell^2(\mathbb N)$ commuting with
    $S+S^*$ and $S+S^*$ would leave invariant a finite dimensional subspace and hence have an eigenvalue. But
    $S+S^*$ is Voiculescu's semicircular element and is well-known not to have an eigenvalue (this follows immediately
    from a direct proof).
    \end{proof}
    
    \begin{corollary}\label{relcomm}Suppose $\delta >1$. For each $k$, $M_0'\cap M_k=P_{0,k}$ (as an algebra).
    \end{corollary}
    \begin{proof} 
    The element $\alpha=$ \vpic{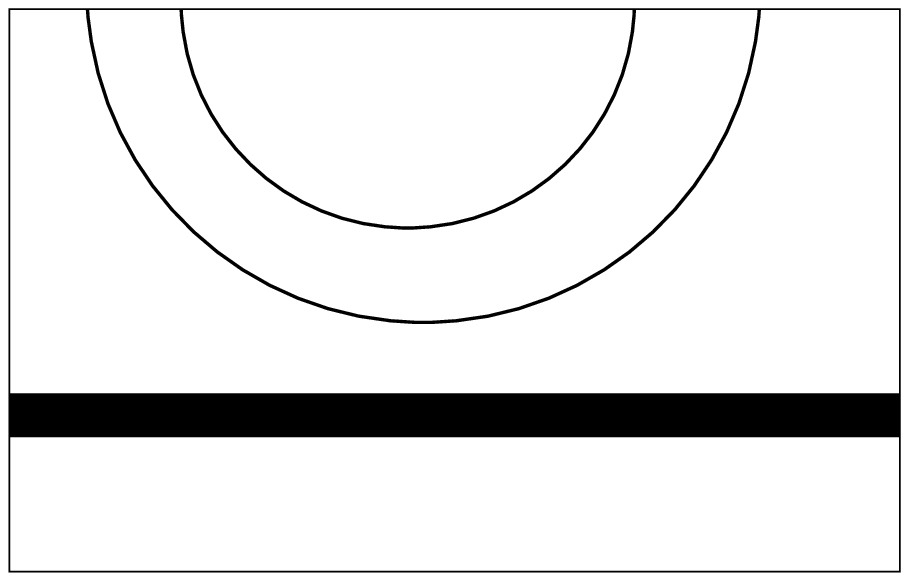} {0.4in}  is in $M_0\subset M_k$
     so it is enough to show that 
    the only elements in the Hilbert space closure of $P_{0,k}A$ that commute with it are elements
    of $P_{0,k}$.  We 
    define  \\$\displaystyle \lambda_n = {1\over (\sqrt \delta)^{n}}  \vpic {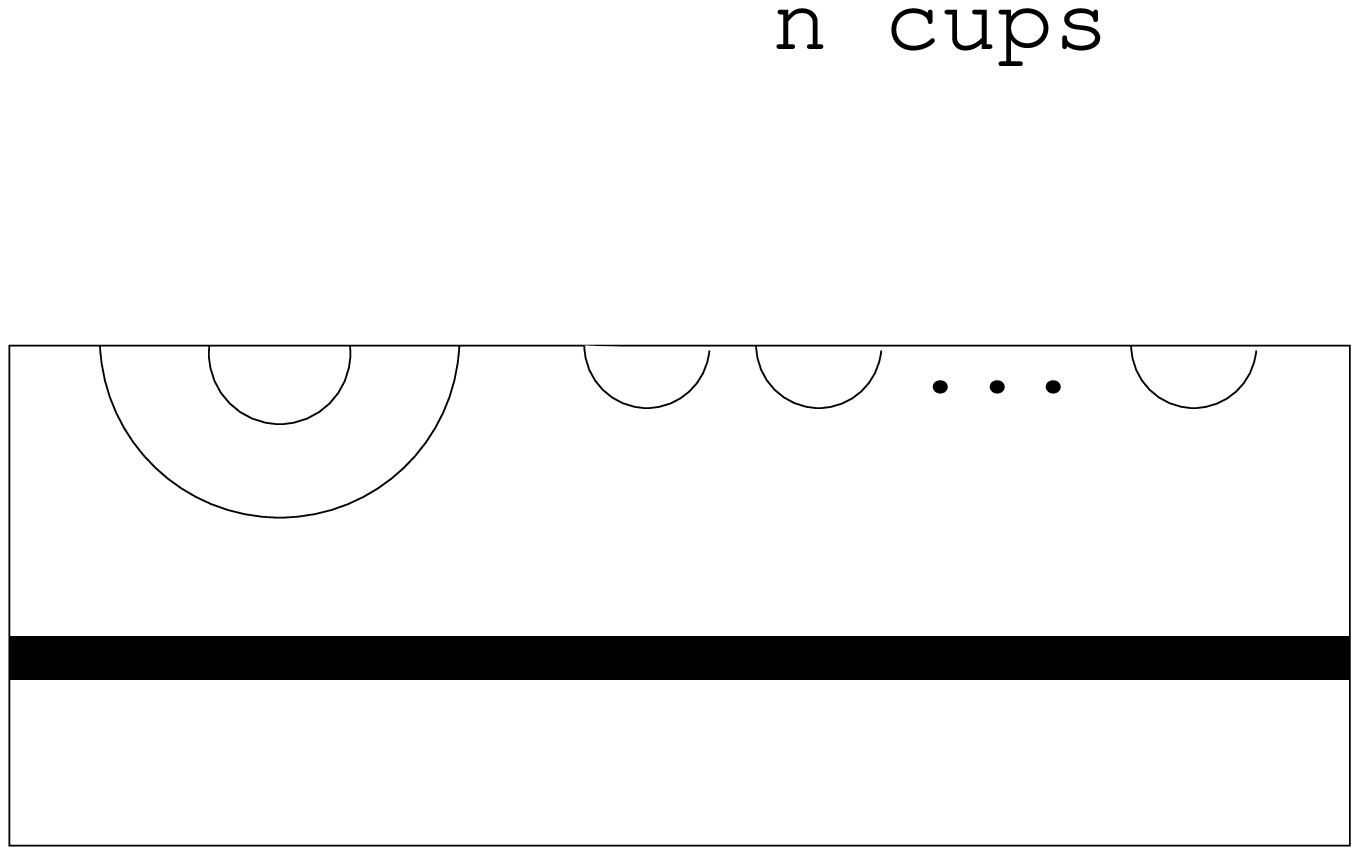} {1.5in} $ 
    and $\displaystyle \rho_n={1\over (\sqrt \delta)^{n}}  \vpic {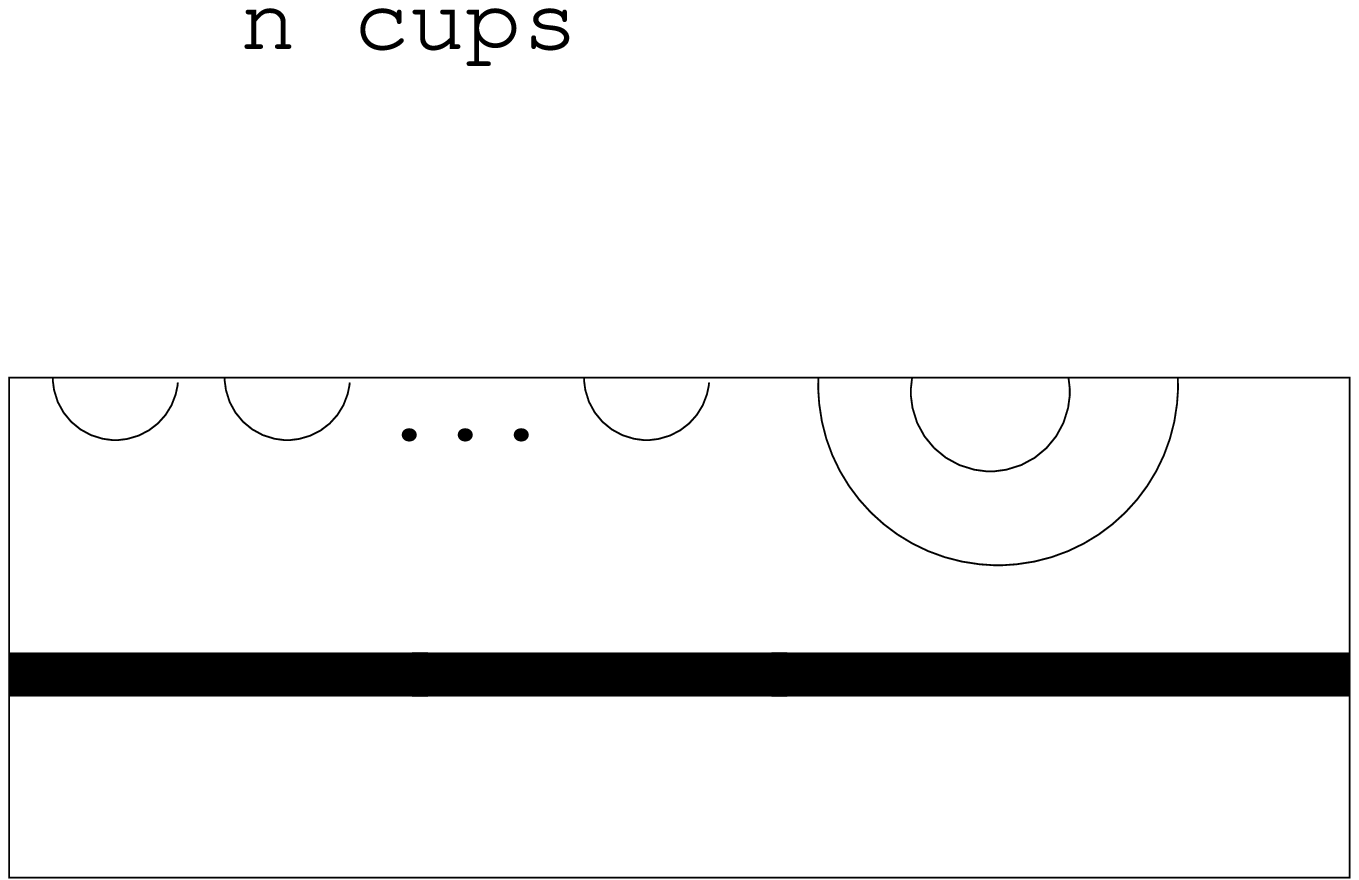} {1.5in} $

    We have
    $$[\alpha,1_{0,n}]=
 (\lambda_n-\rho_n)+\frac{1}{\sqrt \delta}(\lambda_{n-1}-\rho_{n-1})   $$
    
    So if an element $c$ in the closure of $P_{0,k} A$ is written as an $\ell^2$ sum
    $\displaystyle \sum_{n=0}^\infty c_n\star 1_{0,n}$ with $c_n\in P_{0,k}$ we
    find $$[\alpha,c]=\sum_{n=1}^\infty (c_n+\frac{1}{\sqrt \delta} c_{n+1})\star(\lambda_n-\rho_n).$$
    The terms in the sum are orthogonal for different $n$ 
    so for $c$ to commute with $\alpha$ we would have $c_{n+1}=-\delta c_n$ for $n>1$ which forces 
    $c_n=0$ for $n\geq 1$ since  $c_n\in \ell^2$. So $c\in P_{0,k}$.
    \end{proof}
    \begin{corollary} If $\delta >1$, $M_k$ is a type II$_1$ factor.
    \end{corollary}
    \begin{proof} If $x$ were in the centre of $M_k$ it would have to be in $P_{0,k}$. But a trivial diagrammatic argument
    shows that the only elements in $P_{0,k}$ that commute with \vpic {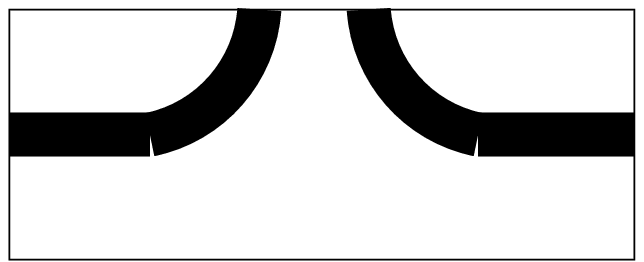} {1in} are scalar multiples of the identity.
    \end{proof}
    
    We now want to identify the II$_1$ factors $M_k$ with the tower coming from the subfactor $M_0\subset M_1$
    obtained by iterating the basic construction of \cite{jones:index}. For simplicity we will do it for the case $M_0, M_1, M_2$, the
    general case being the same argument with heavier notation. 
    \begin{definition}
    The element $e\in M_2$ will be $\frac{1}{\delta}$ \vpic {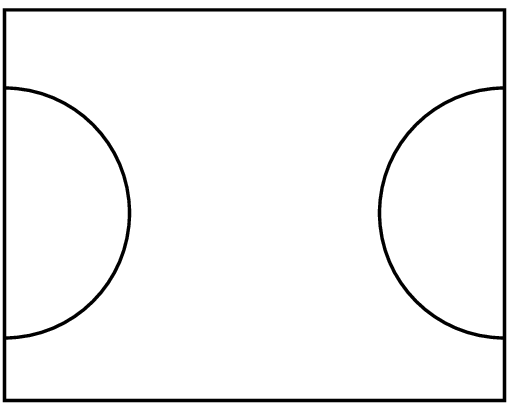} {0.5in} .
    \end{definition}
    \begin{proposition}
    This $e$ is a projection and for $x\in M_1, exe=E_{M_0}(x)e$ and $E_{M_1}(e)=\delta^{-2} id$
    where $E_{M_i}$ is the trace preserving
    conditional expectation onto $M_i$.
    \end{proposition}
    
    \begin{proof} Easy computation with diagrams. \end{proof}
    \begin{lemma}\label{factor} The von Neumann algebra $\{M_1,e\}''$ is  a II$_1$ factor.
    \end{lemma}
    \begin {proof} If $x$ is in the centre of $\{M_1,e\}''$ then it commutes with $M_0$ so by \ref{relcomm} we know
    that $x\in P_{0,2}$. But $x$ also has to commute with \vpic {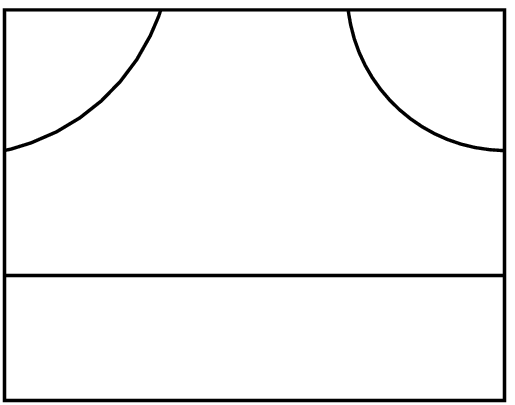} {0.5in} which forces $x$ to be of the from
    \vpic{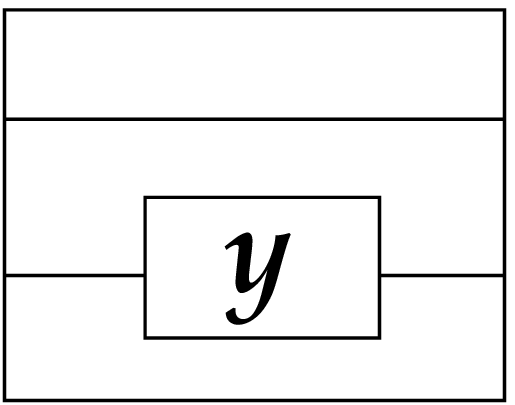} {0.5in} for $y\in P_1$. But for this to commute with $e$ forces it to be a scalar multiple of the
    identity. \end{proof}
   
   \begin{corollary} For $z\in \{M_1,e\}'', ze=\delta^2E_{M_1}(ze)e$.\end{corollary}
   \begin{proof} By algebra, $M$ and $MeM$ span a $^*$-subalgebra of $\{M_1,e\}''$ which is thus
   wealkly dense. The assertion is trivial for $z\in M$ and a simple calculation for   $z\in M_1eM_1$. And $E_{M_1}$ is continuous. \end{proof}
   \begin{corollary} The map $x\mapsto \delta xe$ from $M_1$ to $\{M_1,e\}''e$ is a surjective isometry intertwining
   $E_{M_0}$ on $L^2(M_1)$ and left multiplication by $e$.\end{corollary}
   \begin{proof} Surjectivity follows from the previous lemma. The intertwining property is a calculation.\end{proof}
   \begin{corollary} $\{M_1,e\}''$ is the basic construction for $M_0\subset M_1$ and\\ $[M_1:M_0]=\delta^2$. \end{corollary}
   \begin{proof} The basic construction is the von Neumann algebra on $L^2(M_1)$ generated by $E_{M_0}$ and $M_1$.
   By lemma \ref{factor}, $\{M_1,e\}''$ as a subalgebra of $M_2$ is the same as it is acting on $\{M_1,e\}''e$
   by left multiplication. And this is the basic construction by the previous corollary.
   The index is then just a matter of evaluating the trace of $e$, by uniqueness of the trace on a factor. \end{proof}
    \begin{corollary} $\{M_1,e\}''=M_2$.\end{corollary}
    \begin{proof}  The same argument as above applied to $M_1\subset M_2$ shows that $[M_2:M_1]=\delta^2$. But
    then $[M_2:\{M_1,e\}'']=1$. \end{proof}
    
    Summing up the above arguments applied to the whole tower we have the following:
    \begin{theorem} Let $\mathfrak M_n$ be the II$_1$ factor obtained by the basic construction from 
    $\mathfrak M_{n-2}\subset \mathfrak M_{n-1}$ with $\mathfrak M_0=M_0$ and $\mathfrak M_1=M_1$ and 
    $e_n$ be the projection of the basic construction generating $\mathfrak M_{n+1}$ from $\mathfrak M$.
    Then there is a (unique) isomorphism of towers from $\mathfrak M_n$ to $M_n$
     which is the identity on $M_1$ and sends $e_i$ to $\frac{1}{\delta}$\vpic {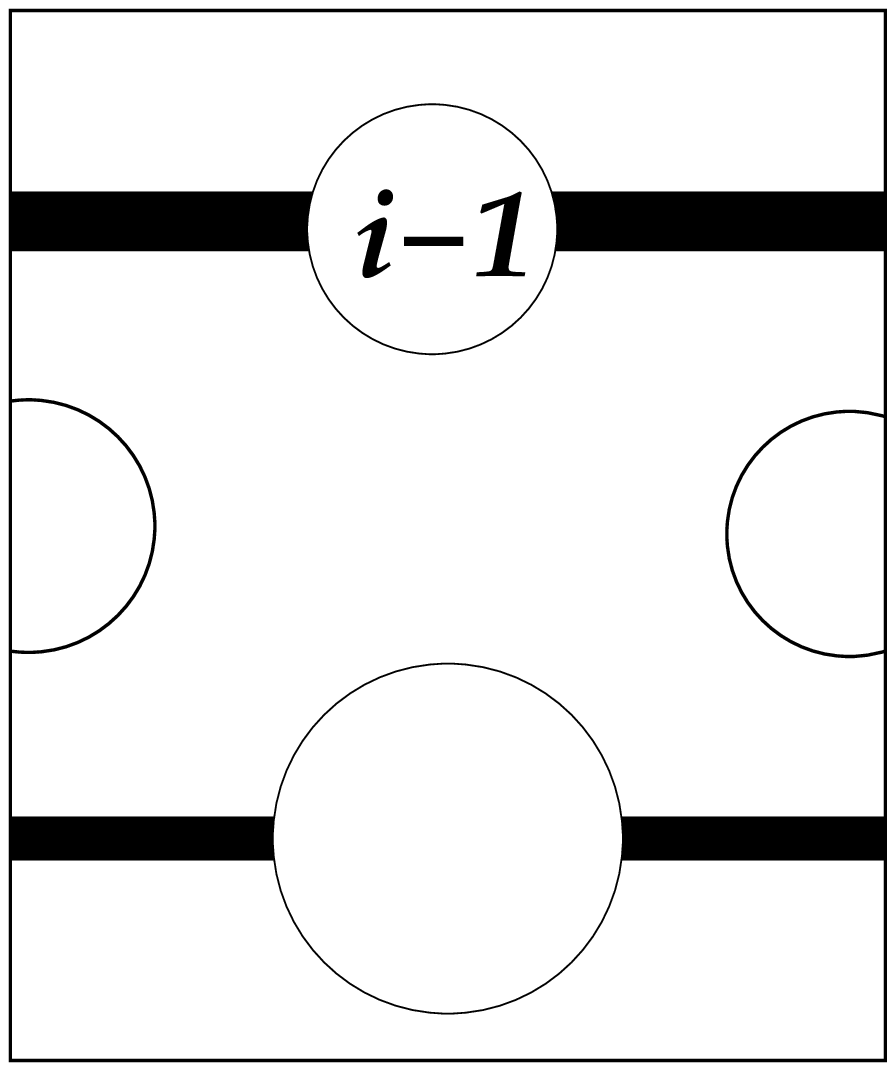} {0.6in} .
     \end{theorem}
     \begin{theorem} Given a subfactor planar algebra $\mathfrak P = (P_n)$ with $\delta >1$ the subfactor $M_0$
     constructed above has planar algebra invariant equal to $\mathfrak P$.
     \end{theorem}
     \begin{proof} It is well known (\cite{jones:planar}) that the planar algebra structure is determined by knowledge of the
     $e_i's$, the multiplication and the embeddings $P_n\subset P_{n+1}$ corresponding to the inclusions
     $M_1'\cap M_{n+1} \subset M_0'\cap M_{n+1}$ and $M_0'\cap M_n\subset M_0'\cap M_{n+1}.$
     But the conditional expectations onto these are just given by the appropriate diagrams. 
     \end{proof}

     \def\nn#1{{{\it \small [#1]}}}
\def\ophs{Hr_k(\mathfrak P)}
\def\fip#1#2{\langle\langle #1, #2 \rangle\rangle}
\def\sip#1#2{\langle #1, #2\rangle}

\section{Change of basis} \label{basis}

In this section we show that the pre-Hilbert space $Gr_k(\mathfrak P)$ defined 
above is isometric and isomorphic as a $^*$-algebra
to the pre-Hilbert space also called 
$Gr_k(\mathfrak P)$
defined in \cite{gjs:free}. To distinguish between them we will call
the latter pre-Hilbert space $\ophs$.

Recall that $\ophs$  is defined
on the same underlying vector space $\oplus_{n\geq 0}P_{n+k}$, but with a simpler multiplication
and more complicated inner product.
The multiplication is simple juxtaposition, $a\bullet b =$ \vpic {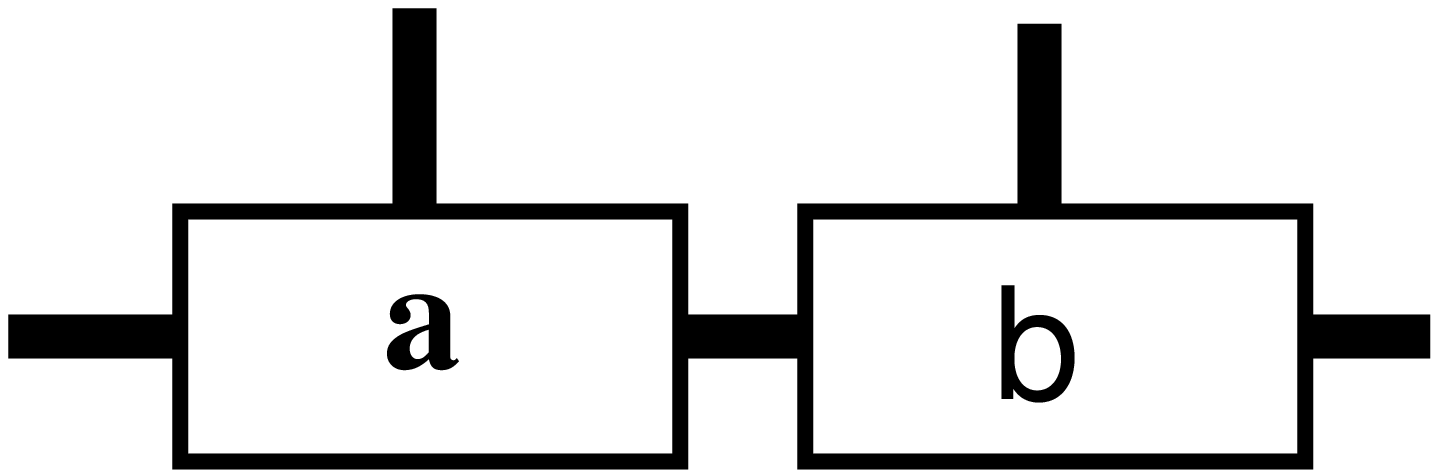} {1in} , while the inner product
$\fip{a}{b}$ of $a \in P_{m,k}$ and $b \in P_{n,k}$ is \vpic {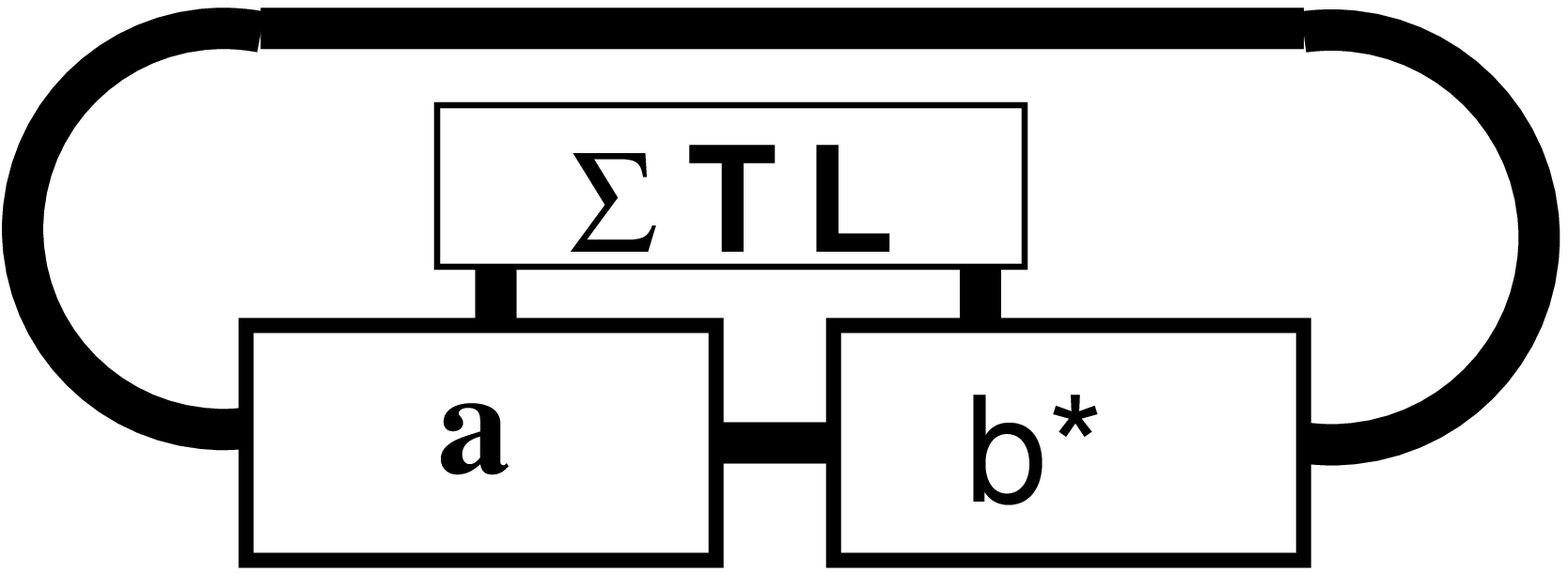} {1in},
where $\Sigma TL$ denotes the sum of all loopless Temperley-Lieb diagrams with
$2(m+n)$ strands on the boundary.
Note that while the multiplication respects the grading, the inner product does not.

We will define an upper-triangular change of basis in $\oplus_{n\geq 0}P_{n+k}$ which
induces an isomorphism between $Gr_k(\mathfrak P)$ and $\ophs$.

Recall that an epi TL diagram is one in which each point on the 
top/outgoing side of the rectangle is connected to the bottom/incoming
side of the rectangle.
A monic diagram is defined similarly, but with the roles of the sides reversed:\\

\noindent \vpic{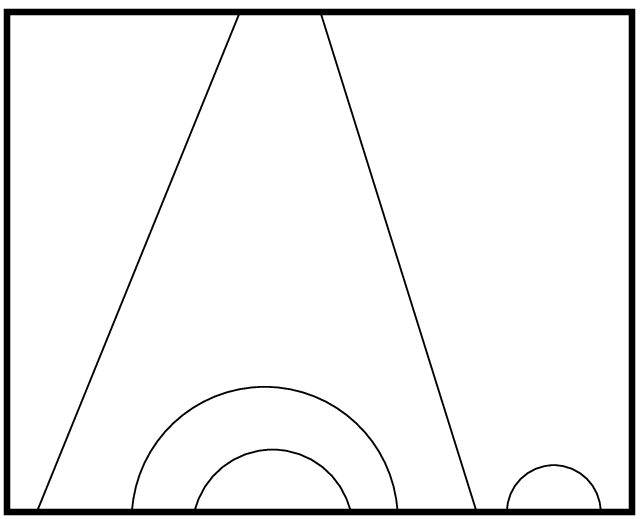} {0.9in}  \vpic{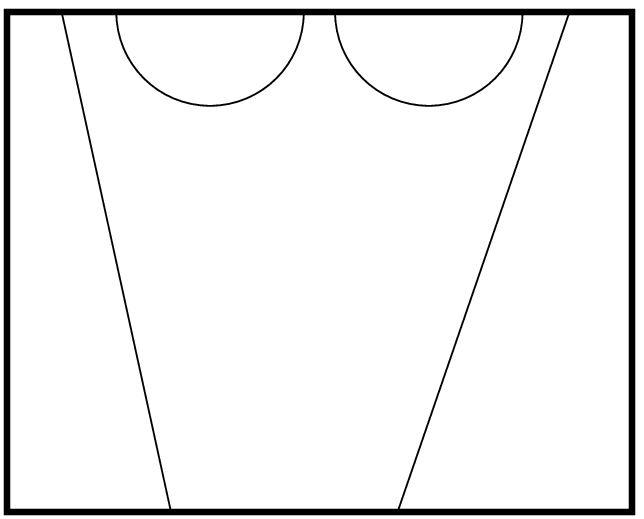}  {0.9in} \hspace{20pt} \vpic{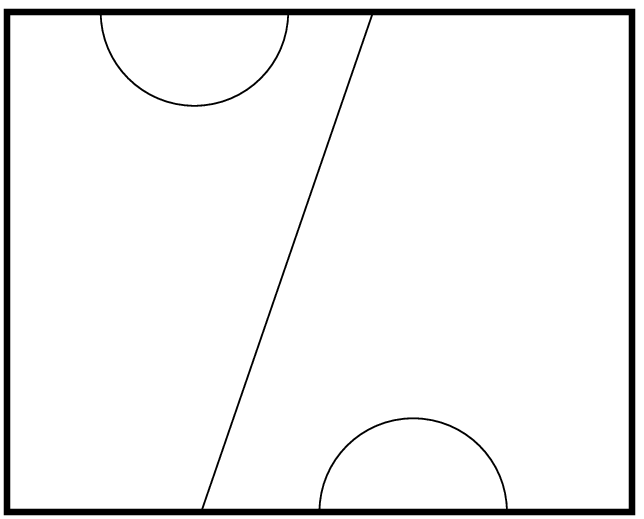} {0.9in} \hspace{45pt} \vpic{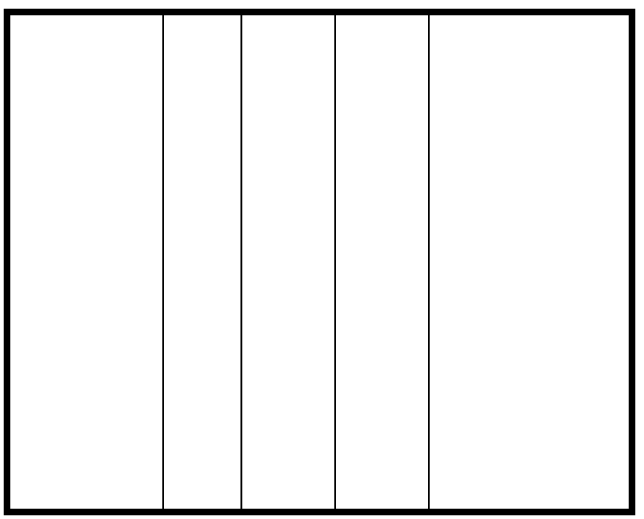} {0.9in} \\
\vspace{10pt}
\qquad epi \hspace{50pt} monic \qquad \hspace{4pt} neither epi nor monic  \qquad {both epi and monic}
  
Note that each TL diagram factors uniquely as epi followed by monic.

We will think of a TL diagram with $2i$ strands on the bottom of the rectangle and
$2j$ strands on the top of the rectangle as a linear map from $P_{i,k}$ to $P_{j,k}$.

Define $X: \oplus_{n\geq 0}P_{n+k} \to \oplus_{n\geq 0}P_{n+k}$ to be the sum of all
epi TL diagrams.
Thus the $j,i$ block of $X$ is the (finite) sum of all epi TL diagrams
from $2i$ strands to $2j$ strands, which is the identity if $i=j$ and zero if $i < j$.

Define a non-nested epi TL diagram to be one where each ``turn-back" or ``cap" on the bottom
of the rectangle encloses no other turn-backs.
Define $Y: \oplus_{n\geq 0}P_{n+k} \to \oplus_{n\geq 0}P_{n+k}$ to be the sum of all
non-nested epi TL diagrams, with the coefficient in the $i,j$ block equal to $(-1)^{i-j}$.

\begin{remark} In the special case of a vertex model planar algebra (\cite{jones:planar}) 
the graded vector space is the (even degree) non-commutative polynomials. Voiculescu
in \cite{voiculescu:symmetries} defined a map from these polynomials to full Fock space, the vacuum component 
of which is the trace on what we have called $\ophs$. in this case the map $X$ gives Voiculescu's
map in its entirety and $Y$ is its inverse. We presume that these formulae are known perhaps in
some slightly different form but have been unable to find them explicitly in the literature.

\end{remark}
\begin{lemma}
$XY = 1 = YX$.
\end{lemma}
\begin{proof}
$X_{jm}Y_{mi}$ is equal to the sum of all products of a non-nested TL diagram from $i$ to $m$
(with $i-m$ turn-backs) followed by a general epi TL diagram from $m$ to $j$,
with sign $(-1)^{i-m}$.
The number of times a given diagram $D$ appears in this sum is is equal to the number
of subsets of size $i-m$ taken from the innermost turn-backs of $D$.
It follows that the total coefficient of $D$ in $\sum_m X_{jm}Y_{mi}$ is
$\sum_p (-1)^p {t \choose p} = 0$ (assuming $p>0$), where $t$ is the total number of 
innermost turn-backs of $D$.
Thus the off-diagonal blocks of XY are zero, and it is easy to see that the diagonal blocks
of $XY$ are all the identity.

The proof that $YX = 1$ is similar, with outermost turn-backs playing the role previously
played by innermost turn-backs.
\end{proof}

\begin{lemma}
$X(a\bullet b) = X(a)\star X(b)$.
\end{lemma}
\begin{proof}
Let $a \in P_{m,k}$ and $b \in P_{n,k}$.
Each epi diagram from $2(m+n)$ to $2j$ appearing in the definition of $X(a\bullet b)$
factors uniquely as $T \cdot (L|R)$, where $L$ is an epi diagram from $2m$ to $2m'$, $R$ is
and epi diagram from $2n$ to $2n'$, $L|R$ denotes $L$ and $R$ placed
side by side, and $T$ is an epi diagram from $2(m'+n')$ to $2j$ where each turn-back has one end
in the $m'$ side and the other end in the $n'$ side.
$L$ corresponds to a diagram used in the definition of $X(a)$,
$R$ corresponds to a diagram used in the definition of $X(b)$,
and $T$ corresponds to a diagram used in the definition of $\star$ 
in $X(a)\star X(b)$.
\end{proof}

\begin{lemma}
$\fip{a}{b} = \sip{X(a)}{X(b)}$.
\end{lemma}
\begin{proof}
Let $a \in P_{m,k}$ and $b \in P_{n,k}$.
Let $D$ be a TL diagram in $\Sigma TL$ used in the definition
of $\fip{a}{b}$.
We can think of $D$ is a TL diagram from $2m$ strands to $2n$ strands, 
and from this point of view it has a unique factorization $E\cdot M$, where $E$ is an epi
diagram starting at $2m$, and $M$ is a monic diagram ending at $2n$.
$E$ is an epi diagram figuring in the definition of $X(a)$,
and $M^*$ is an epi diagram figuring in the definition of $X(b)$.
The way in which $E$ and $M^*$ are glued together corresponds to the definition
of $\sip{\cdot}{\cdot}$.
\end{proof}

\bibliographystyle{alpha}

\end{document}